# Quantitative Genetics and Functional-Structural Plant Growth Models: Simulation of Quantitative Trait Loci Detection for Model Parameters and Application to Potential Yield Optimization.

LETORT Véronique[1*], MAHE Paul[1], COURNEDE Paul-Henry[1], DE REFFYE Philippe[2], COURTOIS Brigitte[3]

[1]Ecole Centrale of Paris, Laboratoire de Mathématiques Appliquées aux Systèmes, F-92295 Châtenay-Malabry cedex, France.

[2]Cirad-Amis, UMR AMAP, TA 40/01 Ave Agropolis, F-34398 Montpellier cedex 5, France and INRIA-Rocquencourt, BP 105, 78153 Le Chesnay cedex, France.

[3]CIRAD-Bios, UMR DAP, Montpellier, F-34398, France.

*Corresponding author: veronique.letort@ecp.fr

## ABSTRACT.

*Background and Aims* Prediction of phenotypic traits from new genotypes under untested environmental conditions is crucial to build simulations of breeding strategies to improve target traits. Although the plant response to environmental stresses is characterized by both architectural and functional plasticity, the recent attempts to integrate biological knowledge into genetics models mainly concern specific physiological processes or crop models without architecture, and thus should prove limited when studying genotype×environment interaction. Consequently, this paper presents a simulation study introducing genetics into a functional-structural growth model, which gives access to more fundamental traits for QTL detection and thus to promising tools for yield optimization.

*Methods* The GreenLab model was selected as a reasonable choice to link growth model parameters to Quantitative Trait Loci (QTL). Virtual genes and virtual chromosomes were defined to build a simple genetic model that drove the settings of the species-specific parameters of the model. The QTL Cartographer software was used to study QTL detection of simulated plant traits. A genetic algorithm was implemented to define the ideotype for yield maximization based on the model parameters and the associated allelic combination.

*Key Results and Conclusions* By keeping the environmental factors constant and using a virtual population with a large number of individuals generated by a Mendelian genetic model, results for an ideal case could be simulated. Virtual QTL detection was compared in the case of phenotypic traits – such as cob weight – and when traits were model parameters, and was found to be more accurate in the latter case. The practical interest of this approach is illustrated by calculating the parameters (and the corresponding genotype) associated with yield optimisation of a GreenLab maize model. The paper discusses the potentials of GreenLab to represent the environment×genotype interactions, in particular through its main state variable, the ratio of biomass supply over demand.

## INTRODUCTION

The main objective of plant genetic studies is to link chromosome loci to specific agricultural traits in the hope of increasing breeding efficiency for crop yield improvement. The recently developed marker-assisted selection strategies rely on attempts to identify and quantify the genetic contributions to the phenotype (set of physical traits). To identify the number and position of loci or genes controlling these target quantitative traits, the overall strategy used by geneticists is to develop a population of individuals (called mapping population) segregating for the target traits and for molecular markers. Markers are "flags" regularly spaced on the whole genome map and representing intergenic (usually non-coding) short strands of DNA that can be hybridised with their counterparts on the target genome, thereby marking a certain location (see Ribaut *et al.*, 2001). Thus it is possible to establish a statistical link between polymorphism at these markers and variability of the target quantitative traits in all individuals of the mapping population. The chromosomal segments, bordered by two adjacent significant markers, are called Quantitative Trait Loci (QTL). They contain the gene of interest but have a confidence interval largely overtaking the gene itself because of the limited power of the classical statistical detection methods. The main phenotypic traits that are classically studied for crops are yield, duration, plant height, resistance to biotic and abiotic stresses, seedling vigour and quality (de Vienne, 1998). Although it has allowed significant advances in crop genetic improvement, there is nowadays a slowdown in yield potential increase for some crops such as rice (Yin *et al.*, 2003). One major difficulty lies in the complex interactions between genotype and environment (G×E) since those traits integrate many physiological and biological phenomena and interactions with field and climatic conditions. Consequently, many QTL are only detected in a narrow range of environmental conditions (Zhou *et al.,* 2007) and the classical genetic models built only from QTL analysis have a correct predictive ability only in a limited range of conditions. It leads to the definition of target environments in breeding programmes and the selection of genotypes adapted to specific environmental characteristics (Hammer *et al.*, 2002). To overcome this difficulty, a growing interest for the use of ecophysiological models is currently emerging; but the communication between those two fields remains difficult. There is an identified need for separating factors influencing a given phenotypic trait and shifting from highly integrated traits to more gene-related traits (Yin *et al.*, 2002). But bridging the gap between genetics models and growth models is still an on-going process, although several studies underlined the potential interest of building such a link (Hammer *et al.*, 2002; Tardieu, 2003; Yin *et al.*, 2004; Hammer *et al.*, 2006).

To deal with the gene level, it seems easier to make the linkage with low-level physiological phenomena. Some attempts to reduce the gap between genetic and ecophysiological models are bottom-up, as in Tomita *et al.* (1999) who simulated the transcription and translation metabolisms for protein synthesis inside a single-cell organism with a virtual genome. But we are still far from getting the whole simulation chain at this level of detail, from the gene expression at the molecular scale to the resulting plant growth processes. Coupland (1995) studied the mutations of the Arabidopsis genome that affect the flowering time and the interactions between genes for the response to long or short days. But he concluded that an accurate modelling was difficult to get because all the genes involved had not been identified yet. Tardieu (2003) argued that using gene regulatory networks to simulate complex gene effects on phenotypic traits was not feasible, due to the large amount of unknown information concerning gene role and regulation rules and to the high number of different genotypes that would have to be analysed.

The top-down approach, considering ecophysiological modelling at a higher organizational level, is more promising. Its principle is to integrate genetic knowledge in plant growth models: for example, Buck-Sorlin (2002) detected QTL for tillering and number of grains per ear in a winter barley population. He used a linear regression to predict the trait values associated with given allelic values at the considered molecular markers and he integrated them into a morphological growth model. But the effect of environment was not taken into account, although it is precisely the role of models to provide helpful tools not only for the dissection of physiological traits into their constitutive components (Yin *et al.*, 2002) but also for unravelling the genotype×environment interactions (Hammer *et al.*, 2005). Dingkuhn *et al.* (2005) tried to link a peach tree model with QTL but the predictive ability of the model decreased when linked with the genetic model. Despite that unconvincing result, their paper illustrates the interest to test further QTL detection for high level model parameters and emphasizes the necessary condition that those parameters should act independently from each other and be subjected to minimal G×E interactions. A successful work was achieved by Reymond *et al.* (2003) who focused on the equation linking leaf elongation rate (LER) to meristem temperature. The three parameters of this equation were fitted from the data and then linked with their associated QTL. Then the link between genetic and ecophysiological models was used to predict leaf elongation rate of non-tested combinations of genotypes and climatic conditions, with satisfactory success (the model explained 74% of the observed variability for LER).

The interest of this approach for breeding strategies is quantified in Hammer *et al.* (2005) using gene-to-phenotype simulations of sorghum: they linked the yield to four basic traits (duration prior to floral initiation, osmotic adjustment, transpiration efficiency, stay-green), the values of which were simulated under three different environmental conditions according to a genetic model built from the relative information found in the bibliography. The simulation results showed that the predictive power and efficiency of marker-assisted selection was enhanced by the link with ecophysiological modelling. They finally discussed the pertinence of such an approach at the plant scale and the level of detail that may be required for the growth model. To add further elements to this discussion, we propose in this paper to examine through a theoretical study the use of a functional-structural growth model as a tool for marker-assisted selection. Since the target traits, such as yield, are the results of the whole plant functioning, it is important to study them in association with all the other processes in the dynamic context of plant growth instead of considering them independently from each other. Functional-structural models (e.g. Wernecke *et al.*, 2000; Drouet and Pagès, 2003; or see Van der Heijden *et al.*, 2007) aim at describing the plant response to environmental factors by integrating ecophysiological functions in the plant architecture at the organ scale. Hence, they can be powerful tools to help analysing the effects of G×E interactions, not pretending that their parameters are directly related to gene expression but assuming that they should, at least, allow detection of more stable QTLs than classically used phenotypic traits. Indeed, parameters for models at organ or plant level already integrate several interacting physiological processes but they are likely to be more stable under diverse environmental conditions than the phenotypic traits that they drive.

Based on this principle, our paper is a first simulation study of QTL detection for parameters of a generic functional-structural growth model on a virtual mapping population built from a simple genetic model. The presentation of this simulation tool of the chain from genotype to phenotype is illustrated with virtual data that allow simplifications to make plant modellers more familiar with the benefits of growth models for breeding work. The formalism of the GreenLab model was chosen: it is a dynamic model taking into account architectural plasticity of the plant and biomass allocation at organ level. Its mathematical formalism allows the easy use of optimization methods, for example in the goal of calculating the best parameters to get an objective criterion under given constraints. This specificity can make it a powerful tool for breeders. The potentialities of such an approach are illustrated using the example of a virtual diploid cereal that could be identified with maize. A genetic algorithm was computed to find the parameters, and therefore the associated genotype, that gives the best yield under a constant environment. This study is the general framework of experiments currently conducted at the Beijing Chinese Academy of Agricultural Sciences (CAAS) on tomato genotypes.

## MATERIALS AND METHODS

### *Main characteristics of GreenLab*

In the model classification proposed by White and Hoogenboom (2003), GreenLab belongs to the level 3 class: it is a generic model whose parameters values are specific for a given plant species. The first step of the work was to add a genetic component: the resulting complete model is thus a class 4 model, *i.e.* a model where "genetic differences are represented by specific alleles, with [allele] action represented through linear effects on model parameters" (White and Hoogenboom, 2003). The resulting flowchart of the final model is represented in Fig. 1. The arrows show the potential influences of the plant genome in the model: it can control the setting of the endogenous parameters of the model and the rules driving the environmental impacts. The circular arrows represent the various feedbacks between organogenesis, biomass production and allocation that can be integrated in the functional-structural model GreenLab. For example, an index of the plant trophic state can drive the architectural development or changes in the plant architecture can induce fluctuations of microclimate.

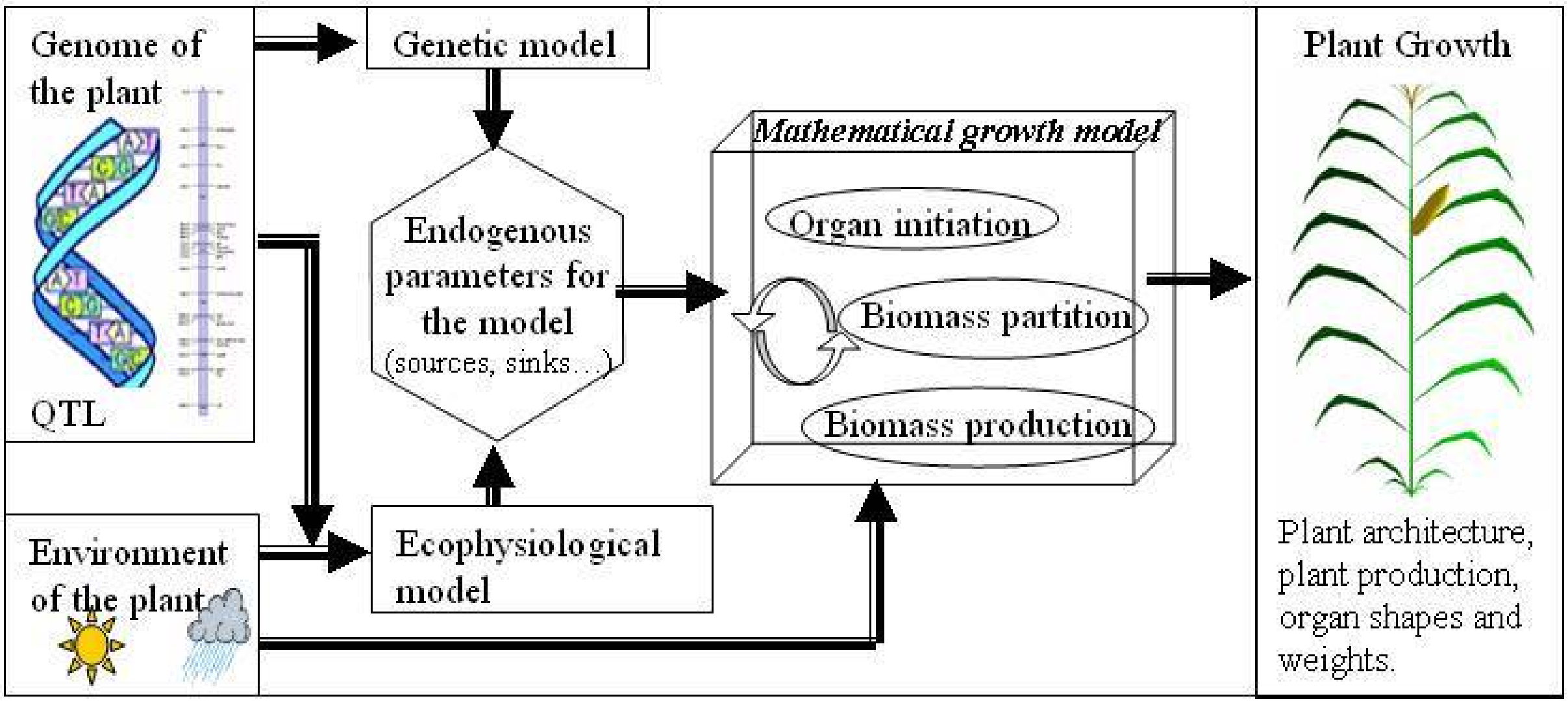

***Figure 1.*** *General flowchart of linking the genetic model to the GreenLab model. The genetic model can potentially have an influence on the determination of species-specific parameters of the model and on the rules driving the environmental impact. The functional-structural model includes complex feedback processes between organogenesis, biomass production and allocation.*

A detailed description of GreenLab can be found in de Reffye *et al.* (1997) and Guo *et al.* (2006). The features useful for the understanding of this study are summarized here.

GreenLab is a generic growth model based on dynamical equations that integrate organogenesis, biomass allocation and production at the organ scale: the plant is regarded as a population of organs classified according to their chronological and physiological ages. Time steps, also called growth cycles, are based on the plant plastochron or phyllochron, *i.e.* are linearly related to the thermal time (Jones, 1992). The net biomass production is computed at each growth cycle and distributed at each cycle to all the expanding organs regardless of their location in the plant (Heuvelink, 1996) and proportionally to their sink strengths. This functional part of the model can be presented in a condensed form through the main recurrence equation that defines the biomass production *Q(n)* of the plant at cycle *n*:

$$Q_n(i) = \frac{E(i)\cdot S_p}{r.k}\left(1-\exp\left(-\frac{P_b}{e}.\frac{k}{S_p}\sum_{i=1}^{t_b} N_b(n-i+1)\sum_{j=1}^{i}\frac{f_b(j)\cdot Q_{n-(i-j)-1}}{D_{n-(i-j)}}\right)\right) \quad (1)$$

This equation is based on the assumption that fresh biomass production is proportional to crop transpiration with an effect of mutual shading of the leaves derived from the Beer-Lambert law (Vose *et al.*, 1995) and adapted to the single plant case (Guo *et al.*, 2006). *E(i)* is the average potential of biomass production during growth cycle *i*, which is determined by the environmental conditions and that can be for example derived from potential evapo-transpiration. The empirical parameter *r* defines a resistance to transpiration, *k* is a factor integrating the light interception effect due to mutual shading of leaves , $S_p$ is the maximal ground projection area available to the plant, $N_b(i)$ is the number of leaves that appeared at cycle *i* and that are photosynthetically active during $t_b$ cycles, *e* is leaf specific weight ($g.cm^{-2}$), $D_i$ is the total plant demand at cycle *i*, which is calculated as the sum of organ sink strengths $P_o$ (in the case of maize, *o* takes its value in the set {b, i, s, c, t} where the letters represent respectively: b: blade, i: internode, s: sheath, c: cob, t: tassel) varying with the organ chronological age according to an empiric function $f_o(j)$ that is defined for each organ type by beta law density function parameters (see Guo *et al.*, 2006 for further details).

The organogenesis simulation relies on the plant decomposition into simple structural units (metamers, axes, structures) and their hierarchical organisation (Barthélémy and Caraglio, 2007). The architecture of the plant is defined by automatic application of rules whose parameters are species-dependent. Those rules can be either predefined (deterministic version of the model, Yan *et al.* (2004)), stochastic (Kang *et al.*, 2003), or dependent on the functional state of the plant (Mathieu *et al.*, 2004). The use of a substructure factorization algorithm also allows a condensed writing of the tree topology by a recurrent procedure (de Reffye *et al.*, 2003; Cournède *et al.,* 2006). Hence, owing to its mathematical formulation, it is possible to study analytically the model behaviour to extract some intrinsic emergent properties (Mathieu, 2006) and to solve optimization problems (Wu *et al.*, 2003). As such, it is a suitable tool for practical applications, such as yield optimization which is one of the main concerns of breeders. This property is illustrated in the following section.

*Genetic model: from genes to model parameters*

This part is a virtual study of the potentials of applying QTL detection methods to GreenLab parameters. To this end, some of the parameters were chosen to be considered as genetically determined and a simple genetic model was built to introduce a plant genotype into the growth model. To illustrate this study, the GreenLab parameters chosen for the simulations are taken from the calibration results of Guo *et al.* (2006) and Ma *et al.* (2007) on Zea *mays* L. The main endogenous parameters can be distinguished on the basis of the stability study made by Ma *et al.* (2007) but here, twelve parameters were arbitrarily chosen:

photosynthetic efficiency, blade thickness, sinks of sheaths, internodes and cob, parameters of sink evolution function for blades, sheaths, internodes and cob, number of shorter internodes at the plant base, cob position on the main stem, seed mass. Those parameters are gathered in an array called *Y*, whose size is *T*, *T* being the number of genetic parameters. Each parameter has a certain range of variation centred on a reference value which is set from the calibration results on maize to obtain simulated phenotypic traits in a valid region.

To simplify the presentation, the virtual genome of the plant is assumed to consist of only one pair of chromosomes, although maize has in reality ten pairs of chromosomes; the general case is easily deduced: the correct chromosome number should be considered if a realistic use of the simulation results was our objective but here, for pedagogic purpose, the clarity of the illustration is privileged. Since all the parameters are quantitative ones, genes can be assumed to be numbers. Each gene can take several values, called alleles. They are written in the matrix *G* (see Fig. 2). The number of alleles for each gene can be easily modified depending on the population studied and is not limited, which allows introducing undetected alleles. Let *N* be the number of genes and *P* the current maximal number of alleles for one gene, then the size of the matrix *G* is ($N \times P$). A chromosome *C* is a vector of size *N* whose components are chosen in the matrix *G* (one allele in each line). The rules driving this choice can be defined by the user in adequacy with the information available about the considered species. For example, it could be necessary to take into account the uneven distribution of genotype frequencies or the skewed distribution of alleles in a natural population, that is generally due to sampling effect due to the small population size, and/or to gametic or zygotic selection in a given area because of the presence of genes influencing gamete or zygote viability, or, in rarer cases, to translocation (e.g. *Musa spec.*, Vilarinhos (2004)). These phenomena could be integrated by setting probabilistic rules to build the individual genotypes. In our paper, the aim is to illustrate QTL detection (the type of mapping population chosen being a recombinant inbred lines population, for which the expected allelic frequency is 1:1 for each individual marker) and potential applications in optimization for selection. Consequently, the alleles are chosen randomly and independently from each other so that all possible genotypes are available. The method is similar to the one adopted by Buck-Sorlin and Bachmann (2000) and Buck-Sorlin *et al.* (2006), except that alleles are considered as variation coefficients (e.g. allele value of 0.9 induces a variation of -10% on the parameters it is related to) instead of integer values. It allows using a simple formalism to define complex rules for the resulting parameter variations.

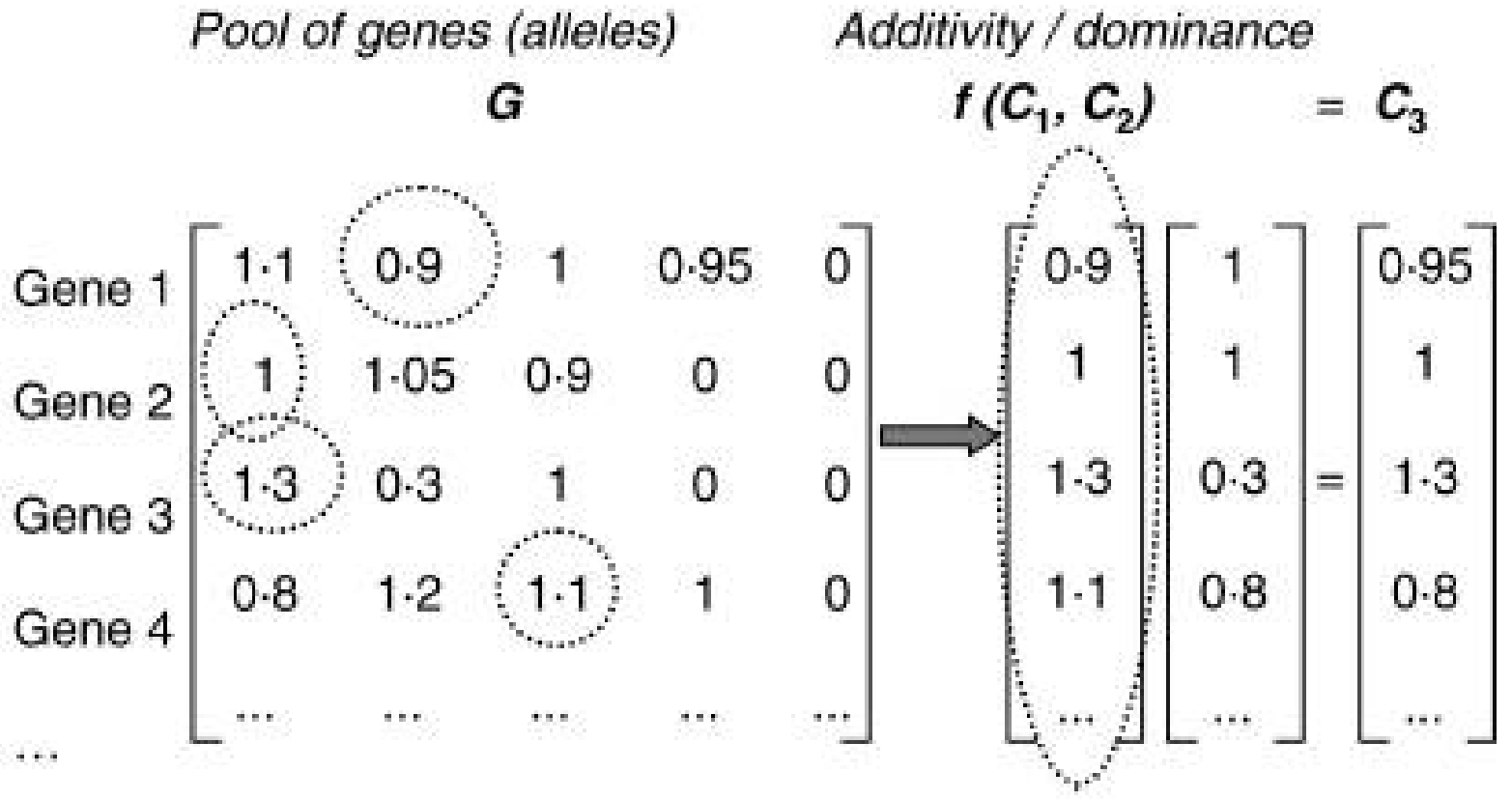

**Figure 2**. *From genes to allele expression. The genotype of the plant is built by choosing alleles among the set of all the possible alleles gathered in matrix G. Function f defines the rules of additivity or dominance that drive the allele effects in the virtual chromosome* $C_3$.

From the values chosen on the pair of chromosomes ($C_1$, $C_2$) of the plant, two kinds of rules drive the effect of those alleles: additivity or dominance. In case of additivity, the resulting effect will be the mean effect of the two alleles, whereas in case of dominance, one allele is chosen to be the one expressed: the choice of the dominant allele is simply represented by their rank in matrix $G$ (the dominant allele is in the first column for each line). The application $f$ is the set of rules for each component of the "chromosome" vectors to get the fictitious chromosome $C_3$ of allele effects, whose size is $N$, by: $C_3 = f\ (C_1,\ C_2)$ (see example in Fig. 2). From that virtual chromosome $C_3$, the 'genetic' vector of parameters is calculated as a product of matrices:

$$Y = D \times A \times C_3 \qquad (2)$$

where $Y$ is the array of the parameters to set and $A$ is a ($T \times N$) matrix defining the influence of genes on each parameter. The matrix $A$ can include pleiotropic rules (one gene has an influence on several parameters) and is also used to define the effect of several genes on one trait (which is the case for quantitative traits). For example, if the first line of matrix $A$ is 2 0 1 0…0, it means that the first parameter depends on the first and third genes; and the influence of the first gene is twice as important as that of the third gene. Epistasis phenomena (effect of one gene on another one) are not considered here. $D$ is a diagonal matrix whose size is ($T \times T$) and whose coefficients are scaling factors to have range compatibility. Indeed, the $j^{th}$ parameter is defined by its variation around its reference value $Y_r(j)$ so the diagonal coefficients $D(i,i)$ of matrix $D$ are defined as :

$$D(j,j) = \frac{Y_r(j)}{\sum_{k=1}^{n} A(j,k)} \qquad (3)$$

The reference value $Y_r$ can be for example the mean value of the parameter in the population.

*Genetic model: simulation of plant reproduction*

The reproduction mechanisms are defined for a diploid plant, that is to say a plant having pairs of homologous chromosomes. For each pair of chromosomes, the 'child' inherits one chromosome from each of its parents. This inherited chromosome can be the result of a crossing-over (exchange of two segments) between the homologous chromosomes of the corresponding parent. Within a population of chromosomes, the number of crossing-over between two markers determines the number of recombinants and is a function of the distance between the two markers. It is assumed here to follow a Poisson law and the points where the cutting occurs are chosen randomly.

The previous section introduced the matrix $A$ that represents the effect of genes on the model parameters. For real experiments, determining the values of the coefficients of matrix A is analogous to QTL detection on model parameters since it relates to searching the associations between locations on genome and parameter values. In our study, the model is used to simulate the phenotypic values and the detection of QTL for the endogenous parameters of GreenLab. For application to real plants, a preliminary step would thus be the estimation of the hidden parameters of the model from the organ- or compartment-level experimental measurements on plants (Ma *et al.*, 2007). Several software packages are used by geneticists to detect QTL, such as QTL Cartographer (Basten *et al.*, 2005). In the simulation, the detection of QTL associated to given traits was done on a mapping population that was generated from recombinant inbred lines: the procedure can be represented as in Fig.

3. First, two individuals are chosen to be the parents, generally with the criterion of being as different and complementary as possible for the considered traits. In the ideal case, those two parents are completely homozygous (*i.e.* same allele values for all genes) so that all individuals issued from their reproduction have the same genome: one chromosome from one parent line (noted 1111…) and one chromosome from the second parent line (noted 2222….). From that $F_1$ generation, several selfings are done until a population whose individuals are homozygous for almost all their genes (97% for the $F_6$ generation) is obtained. To study a real population, the measurements are done on that $F_6$ generation: geneticists genotype each plant with molecular markers covering the whole genome, and measure the target quantitative traits such as ear weight (details can be found in de Vienne, 1998).

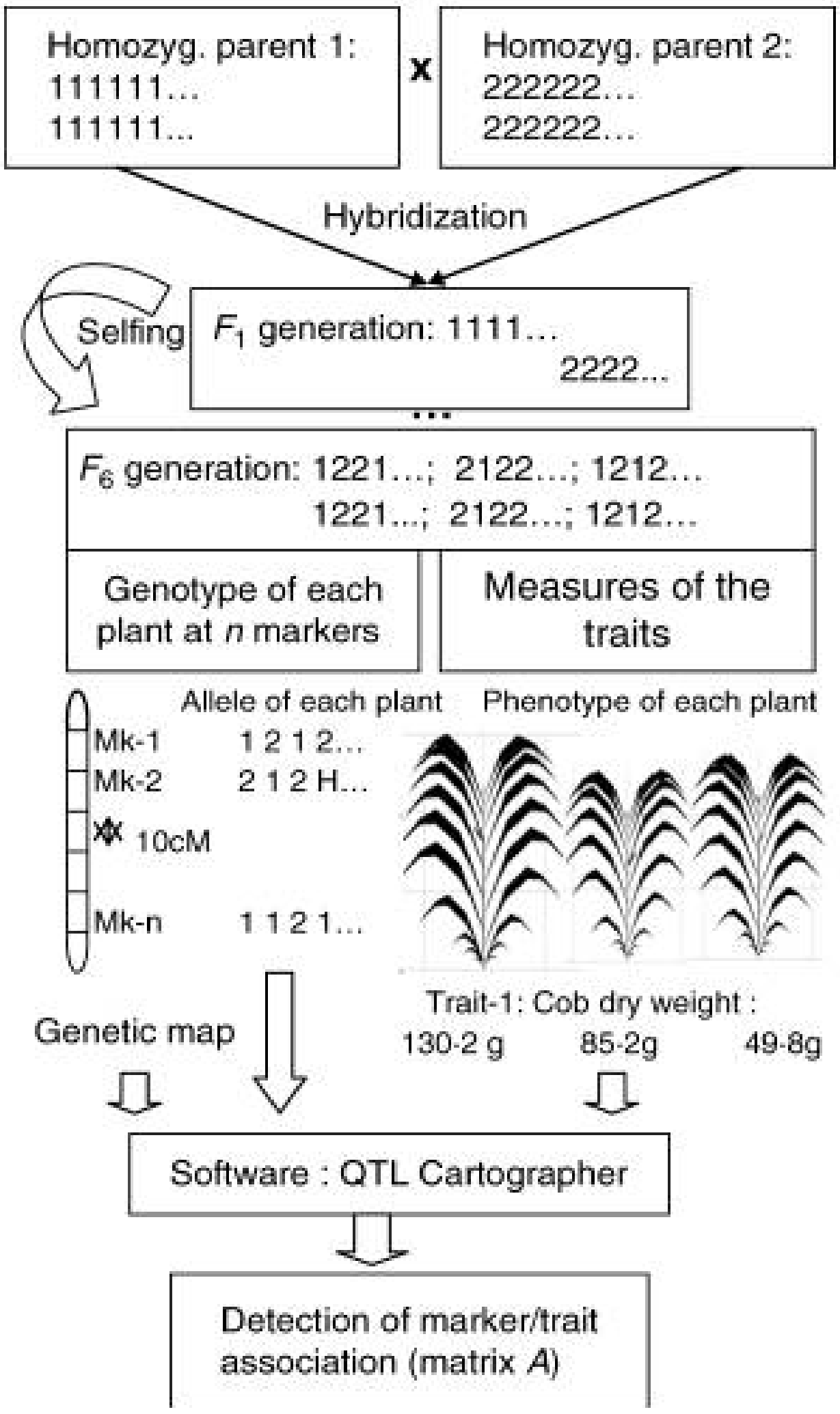

**Figure 3.** *Procedure to build data for QTL detection using QTL Cartographer with recombinant inbred lines. From hybridization of two homozygous parents, an F6 population is obtained by selfings until the sixth generation. Then three kinds of data are collected: molecular map of the genome with distances between markers, genotype of individuals at each marker and phenotype of the same individual for the target trait (e.g. cob weight). From these inputs, single marker analysis is performed using QTL Cartographer to get marker-trait associations.*

For the simulation, the DigiPlant software developed by the Laboratory of Applied Mathematics at the Ecole Centrale Paris (Cournède *et al.,* 2006) was run. The virtual genome of each plant was kept in memory, providing a direct access to its GreenLab parameters and to any phenotypic trait by simulating the plant growth. Thus the three kind of data needed as inputs of QTL Cartographer are gathered for the virtual population (see Fig. 3): *(i)* the genetic map with distances between markers, *(ii)* the genotypes of individuals at all markers (noted 1 for two alleles from parent 1, 2 for two alleles from parent 2 and H for heterozygous marker) and *(iii)* the phenotype of the same individuals for all targeted traits.

To illustrate the potential applications of linking genetic and growth models, a genetic algorithm was computed to find which association of alleles gives a plant with the highest cob weight. The principle of genetic algorithms is derived from the Darwinian rules of genetics of populations: a short introduction can be found in Koza (1995) and Sastry *et al.* (2005). In our study, a simple version was implemented, using the genetic processes defined in the previous section:

(1) An initial population is randomly created, each individual being attributed its genome (which is the "chromosome" vector filled with alleles coding the variables to optimize) and a fitness value (the objective function: the cob weight in our case).

(2) At each iteration, the current population is replaced by a new population, generated with the following steps:

(2.1) Pairs of individuals are selected in the population with a probability depending on their fitness value (this method is called "roulette-wheel selection" in Sastry *et al.* (2005)).

(2.2) These selected individuals can reproduce through crossing-over process ("one-point cross-over"), with a given probability $p_c$.

(2.3) Mutation (change of one allele into another one) can occur with probability $p_m$.

(3) When the final number of iterations is reached, the individual having the best fitness value represents a local solution of the optimization problem.

Thus the average cob weight of the population increases generation by generation thanks to the mechanisms of genetic selection.

## RESULTS

### *QTL detection on GreenLab parameters*

This section presents the results obtained from QTL Cartographer with the set of virtual data, focusing on the comparison of the QTL detection associated to phenotypic traits and to GreenLab parameters.

In this simulation example, the matrix *A* was of size (12×15), *i.e.* $T = 12$ parameters were genetically determined by a set of $N = 15$ genes. Each QTL corresponded exactly to one virtual gene and it was placed at a marker location. Markers were regularly spaced all along the chromosome with a distance of 10 cM between two consecutive markers. Again, this could be changed when considering real mapping data. But nowadays, there are enough markers available in many species (e.g. in Ahn and Tanksley, 1993 or Dunforda *et al.*, 2002) to make a choice of markers regularly spaced. For QTL detection, because of the lack of precision on QTL position (between 10 cM and 30 cM for the QTL confidence interval), one marker every 10 cM is considered enough. To clearly distinguish the QTL, three markers were intercalated between two successive QTL. Single marker analysis was sufficient to detect QTL, since in this virtual study, QTL were represented by the position of non zero components of the matrix *A* defining the influence of genes on the parameters. The single marker analysis method uses linear regression to test the presence of a QTL at each marker by using a likelihood ratio test whose statistic can be converted into a LOD (Logarithm of odds) score as in Eqn (4):

$$LOD = -\log \frac{L_0}{L_1} \qquad (4)$$

where $L_0/L_1$ is the ratio of the likelihood under the null hypothesis (there is no QTL in the interval of markers) to the alternative hypothesis (there is a QTL in the interval). The first trait selected is the first parameter of the model, that is to say the first component of the vector *Y*.

If the first line of the matrix *A* is: (0 0 1 0 0 0 0 1 0 0 0 0 0 0 0), then the LOD curve showing the probability of QTL presence at a marker is presented in Fig. 4A. The position of the two detected QTL is denoted by grey triangles. The LOD scores are very high because, in the genetic model presented in the previous section, alleles have a linear effect on the parameter values. When the trait is a parameter depending on three QTLs with different weights, like in the second line of the matrix *A*: (0 0 3 0 0 0 0 2 0 0 0 0 0 1 0), it gives the curve shown in Fig. 4B. Those examples illustrate that, as expected since virtual data are considered, QTL controlling the endogenous parameters of GreenLab are correctly detected by QTL Cartographer.

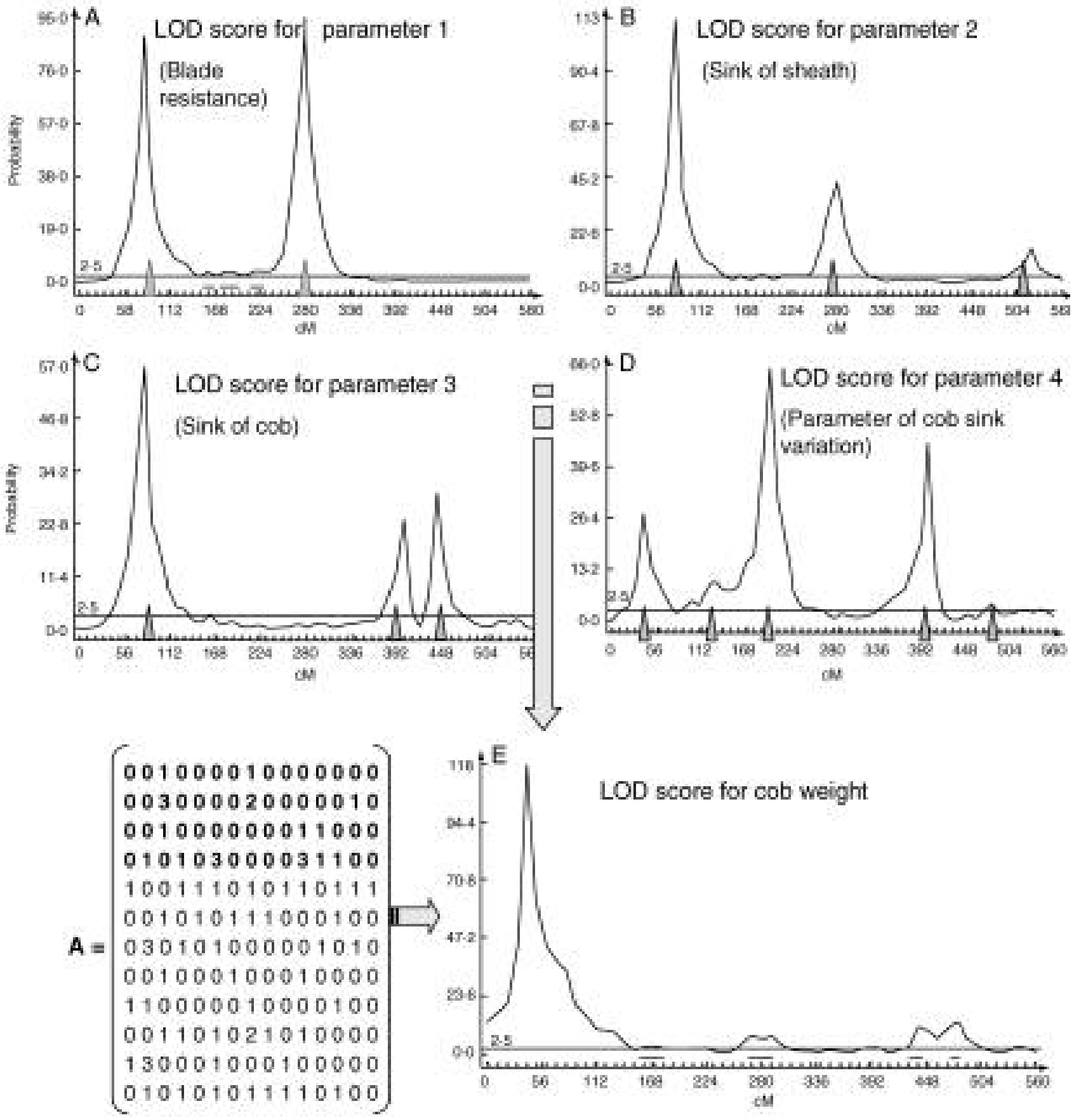

**Figure 4.** *QTL detection on four model parameters (Y(1), Y(2), Y(3), Y(4)) and on the corresponding cob weight. The curves show the probability of QTL presence at each marker position along the chromosome (X-axis represents marker positions in cM). The matrix A coefficients define the effect of each gene on the model parameters. Grey triangles indicate the most probable QTL positions.*

The simulation also allows evaluating what the maximal variation in parameter estimation errors could be that still permits QTL detection. For parameter 8 ($8^{th}$ line of the matrix *A* in Fig.4) related to three major QTL, a coefficient of variation of 15% on the associated values of this parameter decreased sharply the detection, as shown by the comparison of graphs 5A and 5B in Fig.5. This value is in fact a maximal limit since the simulation is done under perfect conditions: no environmental variation, linear effect of genes on model parameters, no epistasis effects and Mendelian segregation.

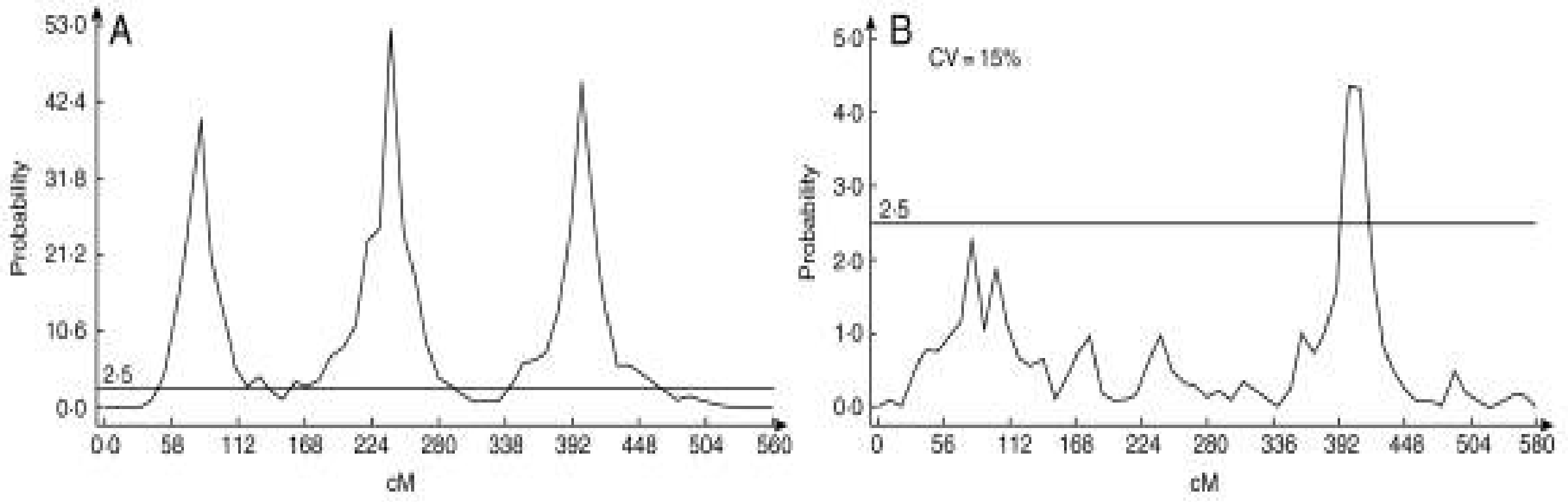

**Figure 5.** *Influence of measurement errors (for parameter values) on QTL detection. The curves show the probability of QTL presence at each marker position along the chromosome (X-axis represents marker positions in cM). With a random white noise on the parameters with a standard deviation equal to 15%, the quality of QTL detection decreased sharply.*

*QTL detection based on phenotypic traits*

The classical direct measurements of plant architecture, got from the growth simulation, were used to feed the process. Cob fresh weight was chosen as a classical phenotypic trait and the relationship between genes and model parameters (matrix *A*) is the one defined in Fig. 4. Fig. 4E gives the results of QTL detection for cob weight: only one major QTL can be detected. The coefficients of the matrix *A* revealed that its position in fact corresponds to genes influencing blade resistance that have a very strong influence on ear weight in the model. However, in graphs 4A, 4B, 4C and 4D other QTL are detected when considering the parameters independently. It means that, for the common measurements done on plant architecture such as plant height, leaf surface or ear weight, only part of the QTL can be detected, even in the ideal case of our simulation. Indeed, those virtual measurements are the result of a step by step plant growth process where all the genetic parameters are involved through complex equations. For example, cob weight at cycle *n* is expressed from Eqn (5) as:

$$W_c(n) = \sum_{i=1}^{n} P_c f_c(i) \cdot \frac{Q_{i-1}}{D_i} \qquad (5)$$

with notations defined in the first section of this paper and the ratio $Q_{i-1}/D_i$ calculated from Eqn (1). It shows that almost all the parameters of the model are involved in the determination of $W_c(n)$. Even under the assumption of constant environment, the values of classical phenotypic traits are the results of complex interacting phenomena that are integrated into the functioning of the growth model. The conclusion of the simulation is that QTL detection gives better results if done on model parameters than on phenotypic traits. Hence growth models can be a useful tool for breeding strategies, but only under the condition that there are ways to control the parameter influence on the phenotypic traits and to optimize their values.

*Determination of the allelic combination optimizing ear weight*

As an example for studying the parameter influence on a phenotypic trait, the relationship between GreenLab parameters and cob weight value was analysed under a stable environment. The coefficients of the matrix *A* are defined by:

$\forall\ i = 1..T,\ \forall\ j = 1..N,\ A(i,j) = 0$ if $i \neq j$ and $A(i,i) = 1$.

It means that each parameter of the model is influenced by only one single QTL. So the detected QTL for cob weight are the ones associated with the parameters that have a strong influence on the calculation of cob weight in the model. Thus it can be seen in Fig. 6 which parameters are the most important for the determination of maize cob weight.

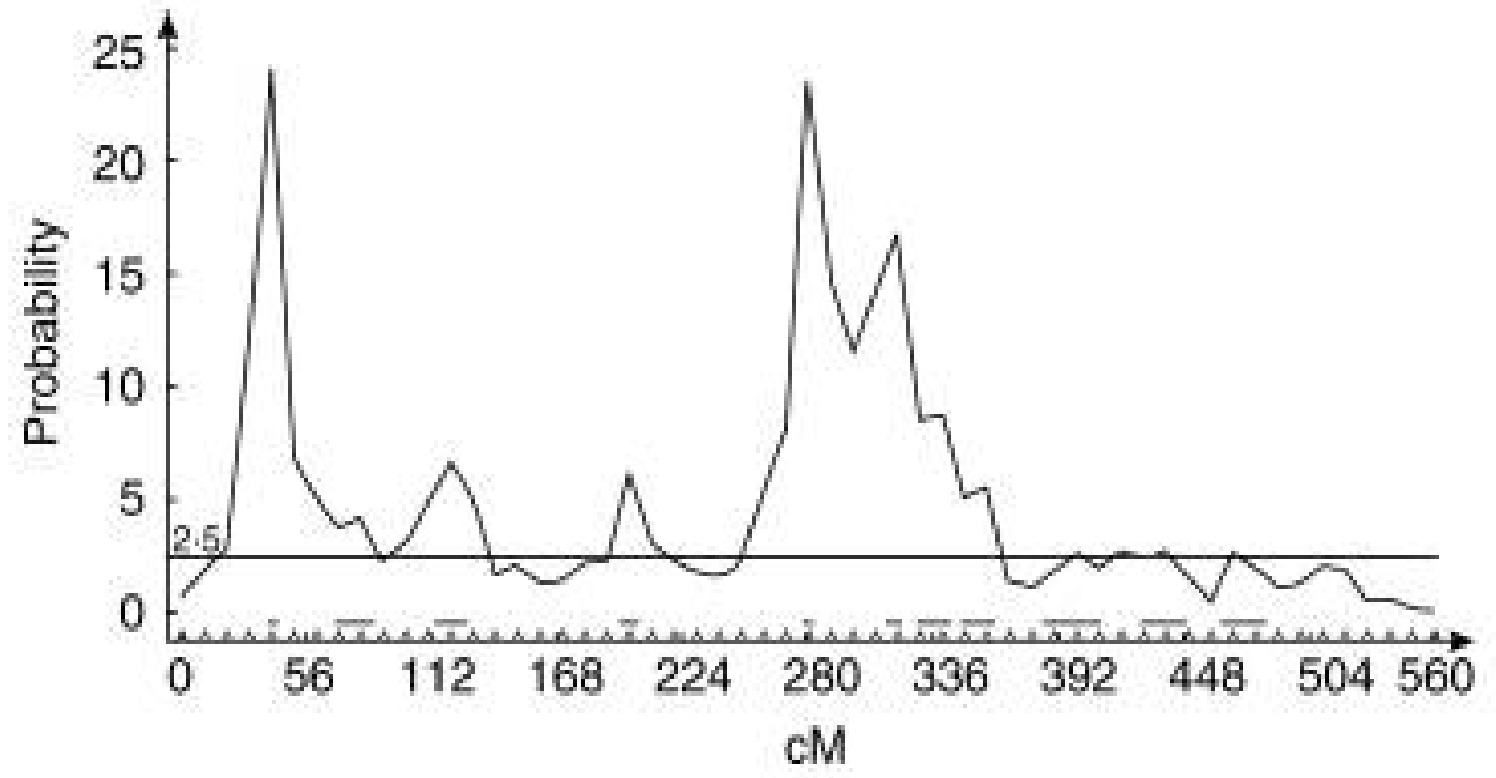

**Figure 6.** *QTL detection for cob weight considering a diagonal matrix A. The curve shows the probability of QTL presence at each marker position along the chromosome (X-axis represents marker positions in cM). Almost all of the 12 genetic parameters of the model are found to have an influence on the cob weight.*

Almost all the QTL positions are detected, which is relevant since all the parameters are linked through Eqn (1) for the determination of cob weight. Moreover, the relationship between cob weight and the model parameters can be complex, as shown in Fig. 7: the shape of the surface defining the cob weight variation is not globally convex.

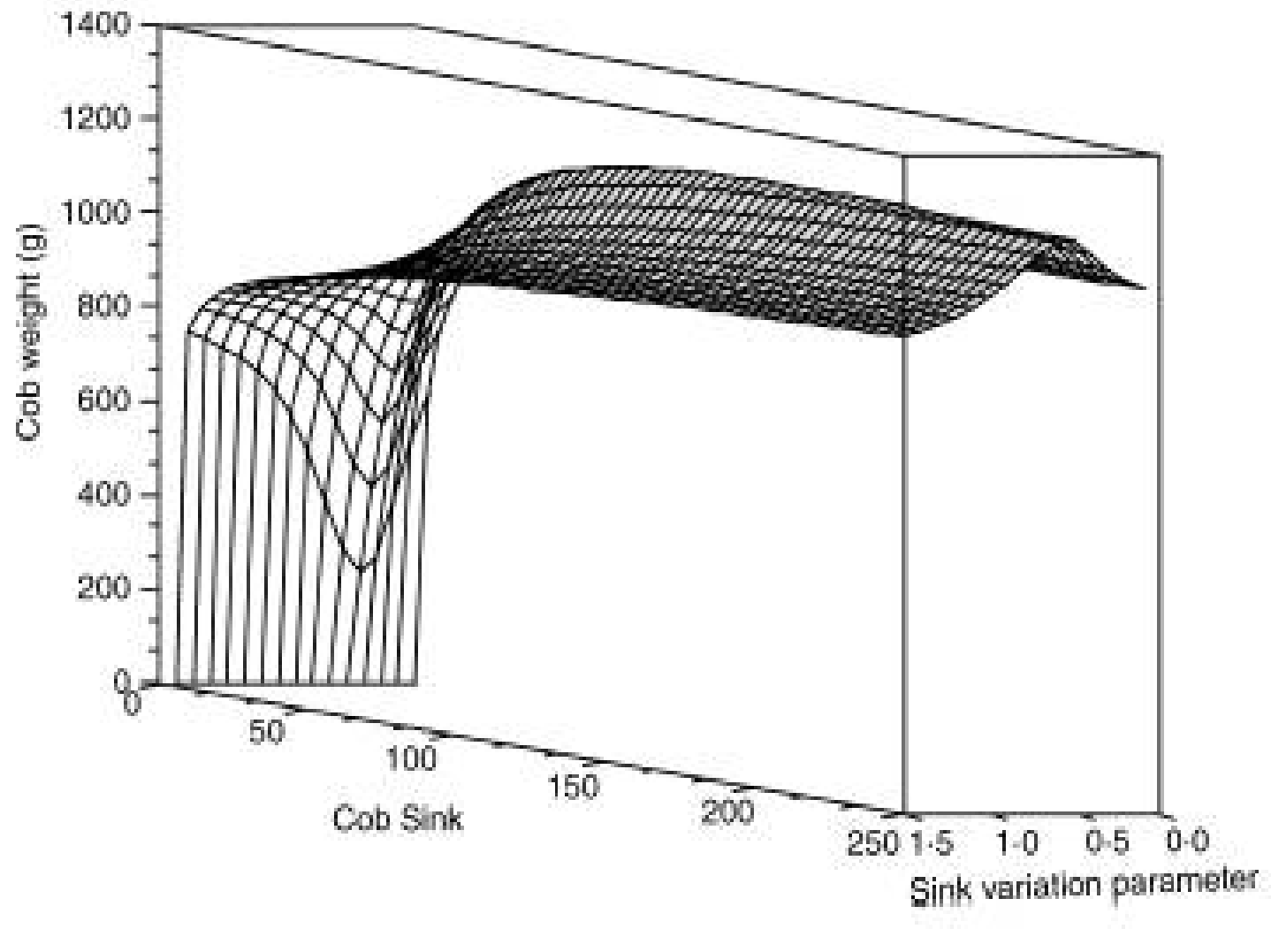

**Figure 7.** *Variation of cob weight according to cob sink and cob sink variation function. Although the surface is not convex, an optimum can be found.*

However, thanks to its mathematical formalism and to its simulation speed, it is possible to apply optimization methods to GreenLab. The results of the genetic algorithm give the allelic combination that optimizes cob weight in given environmental conditions. To simplify the presentation, the matrix *A* was the simplest one: one gene had an influence on only one parameter. It means that the optimization gives the model parameters to get the highest ear weight under constant environmental condition. The procedure could be easily generalized once the coefficients of matrix *A* would have been determined. The twelve parameters that were defined as genetically determined in the first section could take real values, for the ones concerning the functioning part of the model, and integer values for the topological ones. All the other parameters of the model and the environment factors were assumed to be constant. The results are given in Table 1.

**Table 1**. *Optimisation of cob weight under stable environment: parameter ranges and best individual parameter values.*

All parameters are dimensionless except blade thickness (cm) and seed biomass (g). When the optimal value is situated at the interval boundary, it is indicated by "min" (minimal value) or "max" (maximal value). The corresponding cob weight is given from a simulation under the same constant environmental factor.

| Parameter | Reference value | Variation range | Optimal value |
|---|---|---|---|
| Blade thickness (cm) | 0.028 | ±5% | 0.027 *(min)* |
| Blade resistance | 354 | ±5% | 336.3 *(min)* |
| Blade Sink | 1 | - | 1 |
| Sheath Sink | 0.7 | ±10% | 0.63 *(min)* |
| Internode Sink | 2.17 | ±10% | 1.95 *(min)* |
| Cob Sink | 202 | ±30% | 180 |
| Blade Sink Variation Parameter | 0.4 | ±20% | 0.32 *(min)* |
| Sheath Sink Variation Parameter | 0.53 | ±20% | 0.48 |
| Internode Sink Variation Parameter | 0.79 | ±20% | 0.63 *(min)* |
| Cob Sink Variation Parameter | 0.62 | ±30% | 0.7 |
| Number of short internodes at the bottom | 6 | ±20% | 7 *(max)* |
| Cycle of ear appearance | 15 | ±20% | 12 *(min)* |
| Seed Biomass (g) | 0.3 | ±10% | 0.33 *(max)* |
| **Cob weight (g)** | **773** | - | **1221** |

For some parameters, the results could be easily guessed from the analysis of the model behaviour. Blade thickness and blade resistance have to be as small as possible since their diminution increases the plant's ability to perform photosynthesis. On the contrary, large seed biomass gives a stronger plant. Sinks of unproductive organs (except cob) should take minimal values to avoid waste in biomass partitioning. The number of short internodes should be as large as possible since it lets the plant allocate biomass uppermost to the blades that are the future sources of assimilate production. But for other parameters, the influence is more complex and can only be found thanks to the algorithm. The optimization results found for cob sink and the cob sink variation parameters are coherent with the observation of Fig. 7 that tend to show the existence of an optimum point not situated on the interval boundaries. The increase in cob weight induced by the parameter optimization is of about 60%. However, this optimum is not a global maximum, since the use of a genetic algorithm implies that the parameters have only discrete variations in a predetermined space. The grid should be refined and the search domain extended if more precise values were needed for real applications.

## DISCUSSION

In this study, some important aspects of the chain from genetic model to plant growth model was simulated, ending with QTL detection. It is a preliminary step to set the general framework of a simulation tool that will be improved and adapted to specific species when real data of QTL detection on GreenLab parameters are available. The presentation of that simulated procedure could be an original tool to help modellers to understand the potentialities of linking their growth models to quantitative genetics and it illustrates the statement of Dingkuhn *et al.* (2005): "Classical, descriptive phenotyping is based on traits that are too integrative or utilitarian (e.g. yield or leaf area index) and, therefore, insufficiently based on biological functioning to be directly related to gene level information." Indeed, in the simulation, better QTL detection was observed on model parameters than on classical phenotypic traits. Although this study has no real proof value since the simulation was made with hypotheses of simple genetic rules and under constant environment, it has a pedagogic interest since simulation results help to understand the procedure for linking quantitative genetics and ecophysiology and thus enhance communication between those two research fields.

Moreover, it gives the opportunity to discuss further the assets of functional-structural models, and in particular of GreenLab, as candidate plant growth models for QTL detection on their parameters. In the first papers exploring the possibility to link genetic models to plant growth models, the QTL were associated either to the parameters controlling specific physiological phenomena (Reymond *et al.*, 2003; Yin *et al.*, 1999) or to the parameters of crop models (Hammer *et al.*, 2005). However, process-based models present several limitations that could restrict applications in genetics. Indeed, their main drawbacks are: a poor predictive ability of architectural response to environmental factors, such as tillering or organ abortion (Dingkuhn, 1996; Lucquet *et al.*, 2007), difficulties to get reliable computation of leaf area index (LAI) which is mostly the main component of biomass production modules (Marcelis *et al.*, 1998; Heuvelink, 1999), an empirical control of environmental stresses at compartment level (Jeuffroy *et al.*, 2002), difficulties to deal with the inter-plant variability and to handle the often complex interactions between all the different physiological modules (Heuvelink, 1999). These drawbacks result from the fact that process-based models do not take into account plant morphogenesis: at compartment level, since all organs are mixed together, the memory of the growth process is lost and so is the architectural plasticity that reflects the feedbacks between growth and development processes. The endogenous parameters that control both plant development and plant growth are useful key components for yield prediction. Thus they provide new information to renew the breeding process. It provides an adequate strategy to measure plant morphogenesis and to analyze its dynamical biomass production and partitioning.

Several authors (Hammer *et al.*, 2002; Chapman *et al.,* 2003; Tardieu, 2003; Hammer *et*

*al.*, 2006) discussed the properties that growth models should have to expect reasonable chances of success when applied to genetics. Hammer *et al.* (2002) state that their main quality should be a good predictive ability under various environmental conditions. This property can be verified if the growth model parameters define the environmental control of growth phenomena at the different biological levels. Although further analysis still remains to be done, the predictive ability of GreenLab has been demonstrated in Ma *et al.* (2007). The authors found that parameters were stable along development stages and that the model could explain part of the inter-seasonal phenotypic variability. This paper confirmed the analysis of Dingkuhn *et al.* (2005) who discussed the use of GreenLab as a link to genetics. The main drawback they detected was the absence of detailed biological knowledge; however, they suggested that it was "worthwhile to test the GreenLab approach in a genetic context, despite its rudimentary physiology". Indeed, Hammer *et al.* (2002) also emphasized the point that gene-to-phenotype prediction did not require an increase in model complexity, as long as it allowed understanding some key processes so that various combinations of phenotypic responses could be generated through different G×E conditions. The stability analysis of GreenLab parameters tends to reinforce this conviction since it revealed that a small set of chosen rules was sufficient to reproduce plant response to environmental variations (Ma *et al.*, 2007). In the most recent development of GreenLab, it is possible to simulate the complex plasticity of plant architectural and functional responses to environmental factors (Mathieu, 2006). Indeed, the parameters are driven by a state variable of the model: the ratio of global biomass supply $Q$ to total plant demand $D$. The environmental conditions strongly affect the biomass supply and the genetic background of the plant intervenes in the determination of the demand at each growth cycle. That $Q/D$ ratio can be considered as an index of plant vigour and can in particular reflect the environmental impact on plant growth, in interaction with its genome effect. Consequently, the model follows the rules defined by Chapman *et al.* (2003) that stated that a growth model should include "principles of responses and feedbacks" to "handle perturbations to any process and self-correct, as do plants under hormonal control when growing in the field" and to "express complex behaviour (...) even given simple operational rules at a functional crop physiological level".

Another key point is that QTL detection implies heavy data processing on populations of high individual numbers. As in most models, some GreenLab parameters (e.g. organ sinks) cannot be directly measured on plants: those hidden parameters have to be estimated from experimental data collected with destructive measurements. The data collection process for each individual can seem tedious if done on complete measurement (Guo *et al.*, 2006) but, as shown in Ma *et al.* (2007), the number of needed data can be reduced by methods of aggregation or samplings at different levels. Also, the speed of the fitting procedure is a key factor for processing the large size populations required for QTL detection. Thanks to its mathematical formalism, the inverse problem can be computed. GreenLab is associated to a dedicated fitting tool for parameter estimation that relies on the generalized non linear least squares method (Zhan *et al.*, 2003), which allows a very fast resolution (usually, ten iterations are sufficient and the computation time is generally a few seconds).

Finally, it is worthwhile to anticipate what could be the limitations in the use of GreenLab for QTL detection. First, the model's ability to discriminate genotypes with close allelic composition is an important issue (Tardieu, 2003) and depends on the accuracy of the fitting procedure. Also, the level of required accuracy still needs to be determined. Other criteria such as geometrical shape of organs might need to be taken into account, since it is one of the main features used by breeders to differentiate genotypes. In their generic framework for combining crop modelling and QTL mapping to select the best crop ideotype for a specific environment, Yin *et al.* (2003) particularly recommended to test the growth model under several environments: thus the G×E interaction would be analysed in a biological way and not only statistically as in classical genetic models. Concerning the GreenLab model, testing under several environments has been undertaken in Ma *et al.* (2007) but this step should be further investigated.

Moreover, the integrative scale of the growth model may be too large. The basic rules that drive plant growth would thus be unlikely to be the direct expression of independent genes, even if they proved stable in various environmental conditions. Indeed, Luquet *et al.* (2007) investigated the phenotypic impact of a single-gene mutation in the genome of the 'Nipponbare' rice cultivar. They used a model simulating phenotypic plasticity through resource allocation by introducing an internal competition index for the plant. Apart from detailed observations of differences between the growths of mutant and wild cultivars, the estimation of model parameters highlighted that many traits affected by the mutation closely interacted and it was difficult to reconstruct their causal chronology. It means that some traits can be artificially associated to the same QTL even though the underlying gene influences only one physiological function of the plant. Using growth model at the plant level can thus induce artificial pleiotropic effects since the determination of some parameters could be driven by common primary mechanisms (Yin *et al.*, 2003).

A genetic algorithm was used to optimize the parameters in order to get the highest cob weight for maize. One advantage of this kind of optimization algorithm is that it can take into account complex constraints (by defining the viability of individuals) and multi-objective criteria (with weighted fitness values, for example). Thus, if one single allele has combined effects on the phenotype, with positive influence on some traits and negative on others, the algorithm can help to find the best compromise. Here, the optimization procedure was realized on twelve parameters that were considered as genetically determined but in a complete study, more parameters, and their interacting effects, should be included. For example, the importance of tassel presence was not taken into account in the model so tassel sink and its sink variation parameter were kept constant. In the same way, new constraints should be added to have more realistic optimized values. Considering for example plant height, the biomechanical constraints in the internodes were not implemented, thus allometric relationships for internodes were also kept constant and the optimization algorithm gave a sink value for internodes as small as possible. Therefore, the optimization criteria should be adapted and made more complex to answer specific objectives on real species. But it is anyway an interesting contribution of modellers to breeders' work, even if the model relies on simplifying assumptions. The modeller can determine the best allelic combination of genes controlling a given trait through the model under specified conditions. Then the production of the genotype can be more or less difficult depending on the positions of the considered genes and the distances between them, but breeders have developed strategies to separate closely linked genes, involving large segregation populations to get and select the proper recombinant. In any case, it is extremely useful for genotype building to have an idea of the value of virtual ideal genotype without having really to build them, especially in case of pleiotropy when compromises have to be done. This approach could broaden the set of morphological, physiological, biochemical and phenological traits commonly used to characterize ideoptypes, as defined by Donald (1968) and Rasmusson (1987). Using model parameters to build ideoptypes should help overcoming the limitations due to environment pressure on QTL detection (Beattie *et al.*, 2003). Their exploitation in breeding programs, however, is conditioned by their heritability, by the level of genetic variations in the populations and by the genetic correlations among them (Reynolds *et al.*, 2001).

A test of the application of the method is planned to detect QTL for GreenLab parameters on tomato plants. The data collected will feed the simulation tool with real molecular maps, genotypes of individuals and allele effects on the model parameters. A set of experiments is currently done at the Chinese Academy of Agricultural Sciences (CAAS) in Beijing. Tomato plants of about 45 known genotypes are grown in the greenhouse and detailed measurements are done at four growth stages to fit GreenLab parameters. The analysis of those experimental data should provide a further study of QTL detection on model parameters versus phenotypic traits. It should hopefully confirm what this paper only illustrates through simulation, that is to say the potentials of integrating functional-structural models in the gene-to-phenotype chain and the interest of using a mathematical approach to perform optimization processes.

## ACKNOWLEDGEMENT

The authors are grateful to the two anonymous reviewers for their constructive comments on a preliminary version of the paper. This work was supported by the Sino-French Laboratory for Computer Sciences, Automation and Applied Mathematics (LIAMA) in Beijing (China).